\documentclass[11pt]{article}
\usepackage{mathrsfs}
\usepackage[justification=centering]{caption}
\usepackage{booktabs}
\usepackage{multirow}
\usepackage{array}
\usepackage[numbers,sort&compress]{natbib}
\usepackage[colorlinks,
            linkcolor=blue,
            anchorcolor=green,
            citecolor=magenta
            ]{hyperref}
\usepackage{mathrsfs,amssymb,amsfonts}
\usepackage{graphicx}
\usepackage{amsmath,amssymb,latexsym,color}
\usepackage[mathscr]{eucal}
\usepackage[numbers]{natbib}
\newtheorem{thm}{Theorem}[section]
\newtheorem{definition}[thm]{Definition}

\newtheorem{lemma}[thm]{Lemma}

\newtheorem{example}{Example}
\newtheorem{remark}{Remark}[section]

\newcommand{\proof}{{\it Proof.\quad}}
\newcommand{\qed}{\hfill\Box\medskip}
\usepackage{CJK}
\begin{document}
\begin{CJK*}{GBK}{song}
\renewcommand{\abovewithdelims}[2]{
\genfrac{[}{]}{0pt}{}{#1}{#2}}

\title{\bf On $g$-good-neighbor conditional
diagnosability of $(n, k)$-star networks
\thanks{M. Xu's research is supported by the National Natural Science Foundation of China (11571044, 61373021) and the Fundamental Research Funds for the Central Universities.}}

\author{Yulong Wei\quad Min Xu\footnote{\scriptsize Corresponding author.
 {\em E-mail address:} xum@bnu.edu.cn (M. Xu).}\\
{\small  School of Mathematical Sciences, Beijing Normal University,}\\
\noindent {\small Laboratory of Mathematics and Complex Systems,
Ministry of
Education, Beijing, 100875, China}\\ }
 \date{}
 \date{}

 \maketitle

\begin{abstract}
The $g$-good-neighbor conditional diagnosability is a new measure for fault diagnosis of systems. Xu et al. [Theor. Comput. Sci. 659 (2017) 53--63] determined the $g$-good-neighbor conditional diagnosability of $(n, k)$-star networks $S_{n, k}$ (i.e., $t_g(S_{n, k})$) with $1\leq k\leq n-1$ for $1\leq g\leq n-k$ under the PMC model and the MM$^*$ model. In this paper, we determine $t_g(S_{n, k})$ for all the remaining cases with $1\leq k\leq n-1$ for $1\leq g\leq n-1$ under the two models, from which we can obtain the $g$-good-neighbor conditional diagnosability of the star graph obtained by Li et al. [to appear in Theor. Comput. Sci.] for $1\leq g\leq n-2$.


\end{abstract}

\section{Introduction}\label{1}
With the size of multiprocessor systems increasing, processor failure is inevitable. Thus, to evaluate the reliability of multiprocessor systems, fault diagnosability has become an important metric. Many models have been proposed for determining a multiprocessor system's diagnosability. The PMC model was proposed by Preparata, Metze, and Chien \cite{PMC} for fault diagnosis in multiprocessor systems. In the PMC model, all processors in the system under diagnosis can test one another. The MM model, proposed by Maeng and Malek \cite{MM}, assumes that a vertex in the system sends the same task to two of its neighbors and then compares their responses. Sengupta and Dahbura \cite{SD} further suggested a modification of the MM model, called the MM$^*$ model, in which each processor has to test two processors if the processor is adjacent to the latter two processors. Many researchers have applied the PMC model and the MM$^*$ model to identify faults in various topologies (see  \cite{CDH,CH,CHS,CHS16,CLQ,HCS,HL,LL,LTC,LXW,PLT,WHM,WLW,XL,XTH,XTZ,Y,YLM,YLQ,Z,ZQ,ZLX}).

The classical diagnosability for multiprocessor systems assumes that all the neighbors of any processor may fail simultaneously. However, the probability that this event occurs is very small in large-scale multiprocessor systems. In 2005, Lai et al. \cite{LTC} introduced conditional diagnosability under the assumption that all the neighbors of any processor in a multiprocessor system cannot be faulty at the same time. The conditional diagnosability of interconnection networks has been investigated (see \cite{CDH,CHS,CHS16,CLQ,HCS,XTH,XTZ,Y,Z,ZQ,ZLX}).

In 2012, Peng et al. proposed $g$-good-neighbor conditional diagnosability \cite{PLT}, which extended the concept of conditional diagnosability. This requires that every fault-free vertex has at least $g$ fault-free neighbors. Peng et al. \cite{PLT} studied the $g$-good-neighbor conditional diagnosability of the $n$-dimensional hypercube $Q_n$ under the PMC model. Since then, many researchers have studied this topic (see \cite{LL,LXW,WHM,WLW,XL,YLM,YLQ}).

The $(n, k)$-star network $S_{n, k}$, proposed by Chiang and Chen \cite{CC}, is an extension of the $n$-dimensional star graph $S_n$. The network $S_{n, k}$ preserves many ideal properties of $S_n$. In recent years, $S_{n, k}$ has received considerable attention \cite{CDH,CHS16,CLQ,CQS,LGY,LX14,LX,XL,Z}. In particular, Xu et al. \cite{XL} derived the following result about the $g$-good-neighbor conditional diagnosability of $S_{n, k}$ under the PMC model and the MM$^*$ model.
\begin{thm}[Xu et al. \cite{XL}]\label{T1}
The $g$-good-neighbor conditional diagnosabilities of the $(n, k)$-star graph $S_{n, k}$ under the PMC model and the MM$^*$ model are
\begin{equation*}
  t_g(S_{n, k})=
\left\{
  \begin{array}{lll}
   \left\lceil\dfrac{n}{2}\right\rceil-1 & \hbox{if $1\leq g\leq \left\lfloor\dfrac{n}{2}\right\rfloor-1,\ k=1,\ n\geq4$}; \\
   \\
    n+g(k-1)-1 & \hbox{if $1\leq g\leq n-k,\ 2\leq k\leq n-1$},
  \end{array}
\right.
\end{equation*}
and
\begin{equation*}
  t_g(S_{n, k})=
\left\{
  \begin{array}{lll}
  \left\lceil\dfrac{n}{2}\right\rceil-1 & \hbox{if $1\leq g\leq \left\lfloor\dfrac{n}{2}\right\rfloor-1,\ k=1,\ n\geq4$}; \\
  \\
   n+k-2 & \hbox{if $g=1,\ 3\leq k\leq n-1,\ n\geq4$};\\
   \\
       n+g(k-1)-1 & \hbox{if $2\leq g\leq n-k,\ 2\leq k\leq n-1$},
  \end{array}
\right.
\end{equation*}
respectively.
\end{thm}

However, there are some unknown cases (see Table \ref{tab:1}).
\begin{table}[htbp]\scriptsize
\caption{\quad \label{tab:1} The $g$-good-neighbor conditional diagnosability of $S_{n, k}$ under the PMC model and the MM$^*$ model}
\begin{center}
\extrarowheight=4pt
\renewcommand{\arraystretch}{1.3}
\begin{tabular}{|c|c|c|c|c|}\hline
 &$k=1$ &$k=2$ &$3\leq k\leq n-2$& $k=n-1$ \\
\hline
{\multirow{2}{*}{$g=1$}} &$\left\lceil\dfrac{n}{2}\right\rceil-1 ~(n\geq4)$ \cite{XL} & $n$ (PMC) \cite{XL} &{\multirow{2}{*}{$n+k-2$ \cite{XL}}}&{\multirow{2}{*}{$2n-3$~\cite{LL,XL}}}\\

 &$?$ ($n=3$) & $?$ ($n\geq3$, {\rm MM}$^*$)&&\\
\hline
$2\leq g\leq \left\lfloor\dfrac{n}{2}\right\rfloor-1$ &$\left\lceil\dfrac{n}{2}\right\rceil-1$ \cite{XL} &\multicolumn{2}{|c|}{{\multirow{2}{*}{$n+g(k-1)-1$ \cite{XL}}}}&{\multirow{3}{*}{}}\\
\cline{1-2}
$\left\lfloor\dfrac{n}{2}\right\rfloor\leq g\leq n-k$ &$?$ &\multicolumn{2}{|c|}{}&{$(n-g)(g+1)!-1$~\cite{LL}}\\\cline{1-4}

{$n-k\leq g\leq n-2$} &Nonexistence &\multicolumn{2}{|c|}{$?$}&\\
\hline
{ $g=n-1$} &\multicolumn{4}{|c|}{$0$}\\
\hline
\end{tabular}
\end{center}
\end{table}

In this paper, we determine the $g$-good-neighbor conditional diagnosability of $(n, k)$-star networks $S_{n, k}$ for all the remaining cases (see Table \ref{tab:2}). Recently, Li et al. \cite{LL} determined the $g$-good-neighbor conditional diagnosability of the star graph $S_n$ under the PMC model and the MM$^*$ model as follows.
\begin{thm}[Li et al. \cite{LL}]\label{T2}
The $g$-good-neighbor conditional diagnosabilities of the star graph $S_n$ with $n\geq4$ for $0\leq g\leq n-2$ under the PMC model and the MM$^*$ model are $(n-g)(g+1)!-1$.
\end{thm}

Note that $S_{n, n-1}$ is isomorphic to $S_n$. Thus, the results in Table \ref{tab:2} extend their results when $1\leq g\leq n-2$.

\medskip

The rest of this paper is organized as follows. Section \ref{3} introduces some terminology and preliminaries. Our main results are given in Section \ref{4}. Finally, Section \ref{6} concludes the paper.
\section{Terminology and preliminaries}\label{3}
An undirected simple graph $G =\big(V(G), E(G)\big)$ is used to represent a system (or a network) where each vertex represents a processor and each edge represents a link. A {\em subgraph} $H$ of $G$ is a graph with $V(H)\subseteq V(G)$, $E(H) \subseteq E(G)$, and the endpoints of every edge in $E(H)$ belonging to $V(H)$. For an arbitrary subset $F\subseteq V(G)$, we use $G-F$ to denote the graph obtained by removing all the vertices in $F$ from $G$. Given a nonempty vertex subset $V'$ of $V(G)$, the {\em induced subgraph} by $V'$ in $G$, denoted by $G[V']$, is a graph in which the vertex set is $V'$ and the edge set is the set of all the edges of $G$ with both endpoints in $V'$. For a given vertex $v$, we define the {\em neighborhood} $N_G(v)$ of $v$ in $G$ to be the set of vertices adjacent to $v$. The {\em degree} of vertex $v$, denoted by $d_G(v)$, is the number of vertices in $N_G(v)$. The {\em minimum degree} of a graph $G$, denoted by $\delta(G)$, is $\min_{v\in V(G)} d_G(v)$. A graph $G$ is $k$-{\em regular} if $d_G(v)=k$ for any $v\in V$. For a given set $A\subseteq G$, we denote by $N_G(A)$ the set $\big(\bigcup _{v\in V(A)}N_G(v)\big)-V(A)$. For neighborhoods and degrees, we omit the subscripts of the graphs when no confusion arises. The {\em symmetric difference} of two sets $F_1$ and $F_2$ is defined as the set $F_1\bigtriangleup F_2=(F_1-F_2)\cup (F_2-F_1)$. Please refer to \cite{BMG} for graph-theoretical terminology and notation undefined here.

Now we focus on star graphs $S_n$ and $(n, k)$-star networks $S_{n, k}$. For a given integer $n$ with $n\geq 1$, we set $I_n=\{1,2,\ldots, n\}$ and $I'_n=I_n-\{1\}$. Let $P(n, k)=\{p_1p_2\ldots p_k\mid p_i\in I_n, p_i\neq p_j, 1\leq i\neq j\leq k\}$, the set of $k$-arrangements on $I_n$, where $k\in I_n$.  We will abbreviate $P(n,n)$ as $P(n)$.

\begin{figure}[hptb]
  \centering
  \includegraphics[width=10cm]{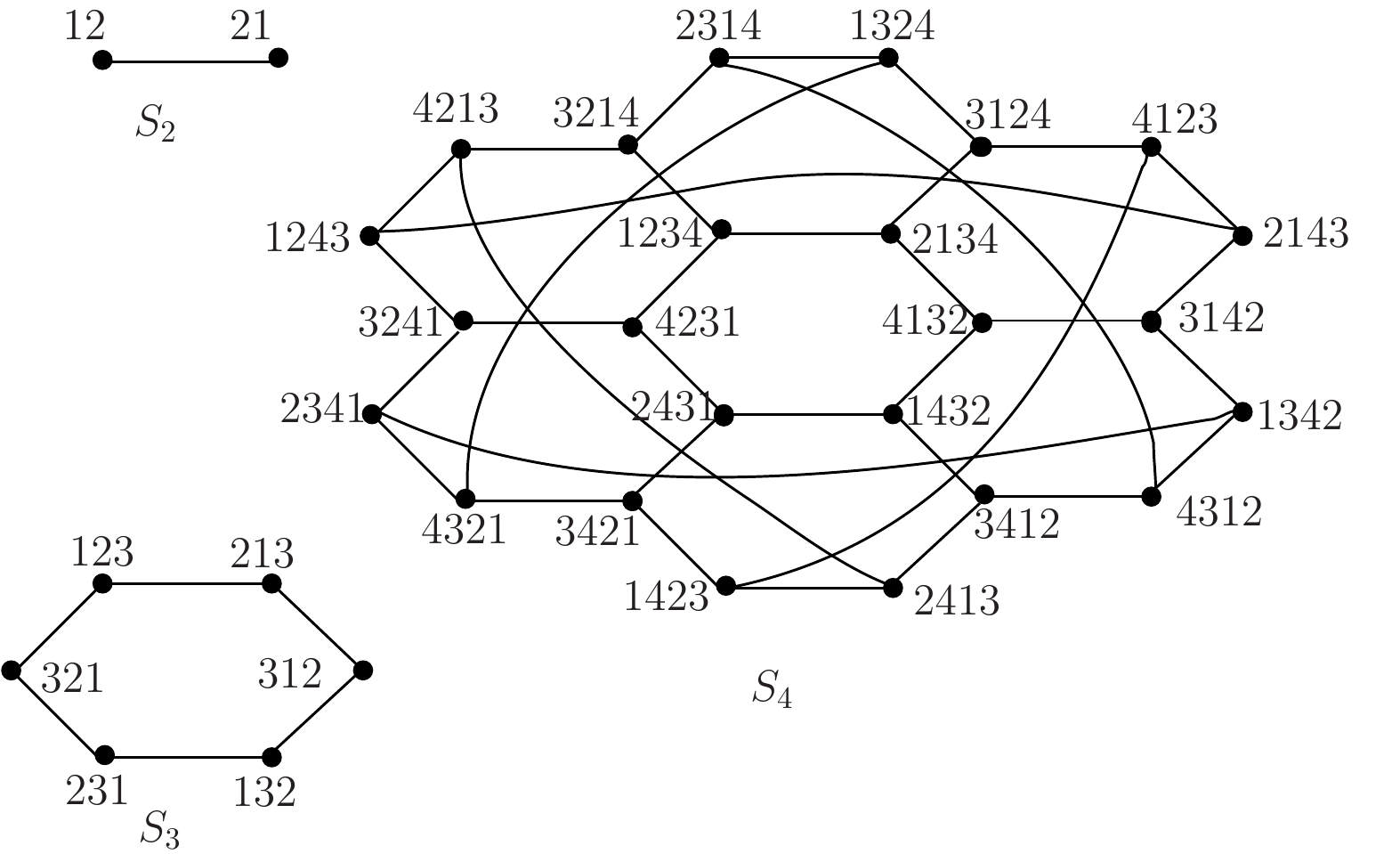}\\
  \caption{Star graphs $S_2$, $S_3$, and $S_4$.
}\label{D001}
\end{figure}

\begin{definition}[Akers and Krishnamurthy \cite{AK}]\label{D01}
An $n$-dimensional star graph $S_n$ is a graph with vertex set $P(n)$, a vertex $p=p_1p_2\ldots p_i\ldots p_n$ being linked a vertex $q$ if and only if $q=p_ip_2\ldots p_{i-1}p_1p_{i+1}\ldots p_n$ for some $i\in I'_n$ (See Figure \ref{D001}).
\end{definition}

\begin{figure}[hptb]
  \centering
  \includegraphics[width=10cm]{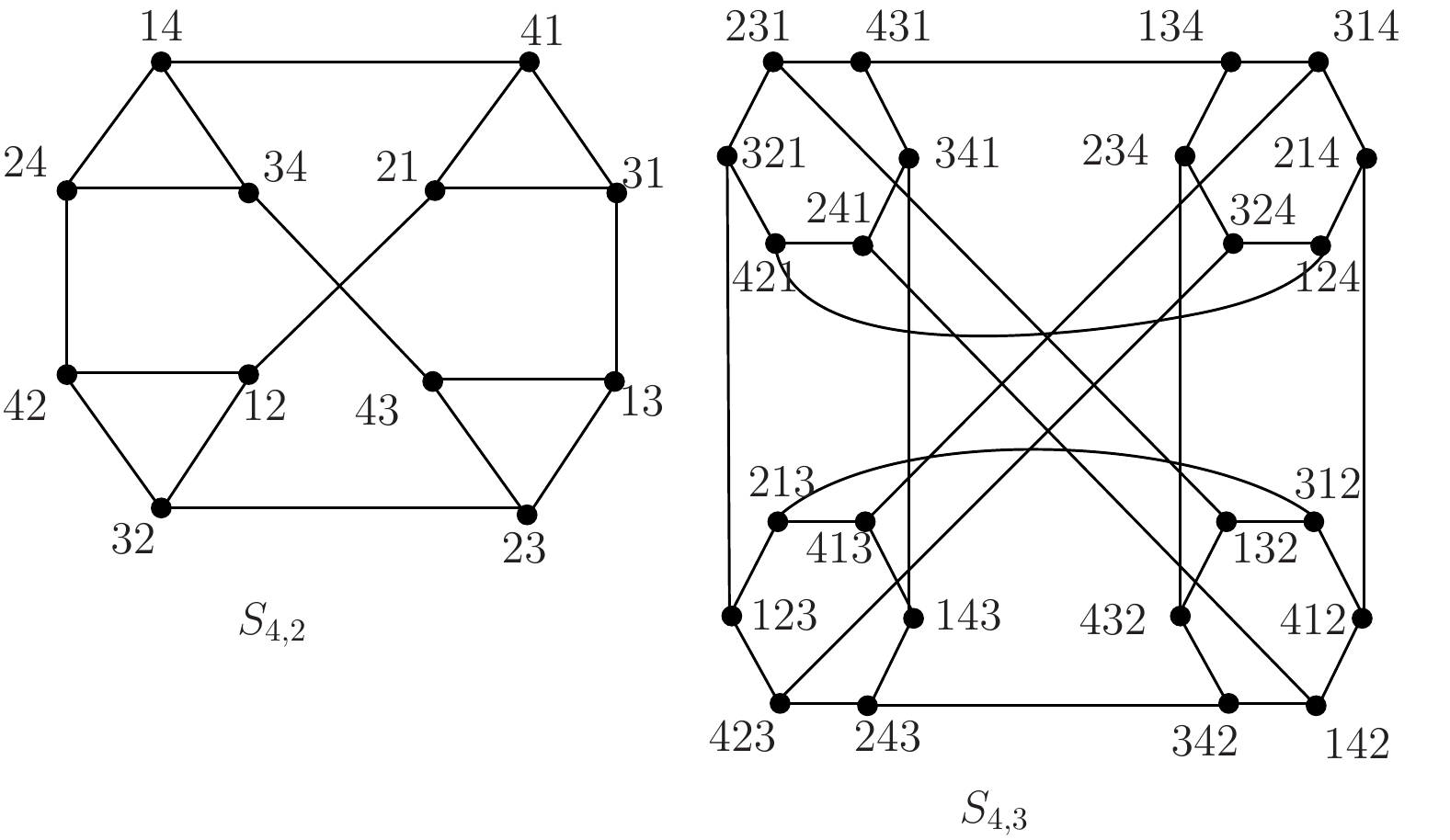}\\
  \caption{$(n, k)$-star graphs $S_{4, 2}$ and $S_{4, 3}$.
}\label{D002}
\end{figure}

\begin{definition}[Chiang and Chen \cite{CC}]\label{D02}
An $(n, k)$-star graph $S_{n, k}$ (See Figure \ref{D002}) is a graph with vertex set $P(n, k)$, a vertex $p=p_1p_2\ldots p_i\ldots p_k$ being linked a vertex $q$ if and only if $q$ is:
\begin{enumerate}
\item[{\rm (a)}] $p_ip_2\ldots p_{i-1}p_1p_{i+1}\ldots p_k$, where $i\in I'_k$ (swap $p_1$ with $p_i$); or

\item[{\rm (b)}] $p'_1p_2p_3\ldots p_k$, where $p'_1\in I_n-\{p_i\mid i\in I_k\}$ (replace $p_1$ by $p'_1$).
\end{enumerate}
\end{definition}

By the definition, the $(n, k)$-star graph is an $(n-1)$-regular, $(n-1)$-connected graph and $S_{n, n-1}$ is isomorphic to $S_n$ (see \cite{CC}).

Now we introduce two models for fault diagnosis.

In the PMC model, all processors in the system under diagnosis can test one another. The set of tests can be represented by a directed graph $G = (V, E)$, in which each vertex represents a processor, and an edge $(u, v)$ indicates that the processor $u$ has tested processor $v$. The outcome of vertex $u$ testing vertex $v$ is denoted by $\sigma(u, v)$, where
\begin{equation*}
\sigma(u, v)=
\left\{
  \begin{array}{lll}
   0 & \hbox{if $\{u, v\}\cap F=\emptyset$}; \\
   \\
   1 & \hbox{if $u\notin F$, $v\in F$}; \\
    \\
    0~{\rm or}~1 & \hbox{if $u\in F$,}
  \end{array}
\right.
\end{equation*}
where $F$ is the set of faulty vertices.

In the MM$^*$ model, a processor executes comparisons for any pair of its neighboring processors. A graph $G =(V, E)$ is used to represent a system, where each vertex represents a processor and each edge represents a link. Assign a task to each vertex. The vertex $w$ is a comparator of a pair of vertices $\{u, v\}$ if $(u, w)\in E$ and $(v, w)\in E$. The outcome of this comparison is denoted by $\sigma\big((u, v)_w\big)$, where
\begin{equation*}
\sigma\big((u, v)_w\big)=
\left\{
  \begin{array}{lll}
   0 & \hbox{if $\{u, v, w\}\cap F=\emptyset$}; \\
   \\
   1 & \hbox{if $w\notin F$, $\{u, v\}\cap F\neq\emptyset$}; \\
    \\
    0~{\rm or}~1 & \hbox{if $w\in F$,}
  \end{array}
\right.
\end{equation*}
where $F$ is the set of faulty vertices.

The collection of all outcomes is called a {\em syndrome} $\sigma$. The diagnosis problem involves using the syndrome to determine the status (faulty or fault free) of each processor in the system. For a given syndrome $\sigma$, a subset $F\subseteq V$ is said to {\em be consistent with} $\sigma$ if the syndrome $\sigma$ can be produced from the faulty set $F$. In the PMC model, $F$ is said to be consistent with $\sigma$ if the syndrome $\sigma$ can be produced from the situation that, for any $(u, v) \in E$ such that $u \notin F$, $\sigma(u, v)=1$ if and only if $v \in F$. In the MM$^*$ model, $F$ is said to be consistent with $\sigma$ if the syndrome $\sigma$ can be produced from the situation that, for any $(u, w)\in E$ and $(v, w)\in E$ such that $w \notin F$, $\sigma\big((u, v)_w\big)=1$ if and only if $\{u, v\}\cap F\neq\emptyset$. Therefore, on the one hand, a faulty set $F$ may produce a number of different syndromes. On the other hand, different faulty sets may produce the same syndrome. Define $\sigma (F)=\{\sigma \mid F ~{\rm is~ consistent~ with}~ \sigma\}$. Two distinct sets $F_1$, $F_2 \subseteq V$ are said to be indistinguishable if $\sigma(F_1)\cap \sigma(F_2)\neq\emptyset$; otherwise, $F_1$ and $F_2$ are said to be distinguishable. We say that $(F_1, F_2)$ is an indistinguishable pair if $\sigma(F_1)\cap \sigma(F_2)\neq\emptyset$; otherwise, $(F_1, F_2)$ is a distinguishable pair.

The following lemmas give necessary and sufficient conditions for a pair of sets to be distinguishable under the PMC model and the MM$^*$ model.

\begin{figure}[hptb]
 \centering
  \includegraphics[width=8cm]{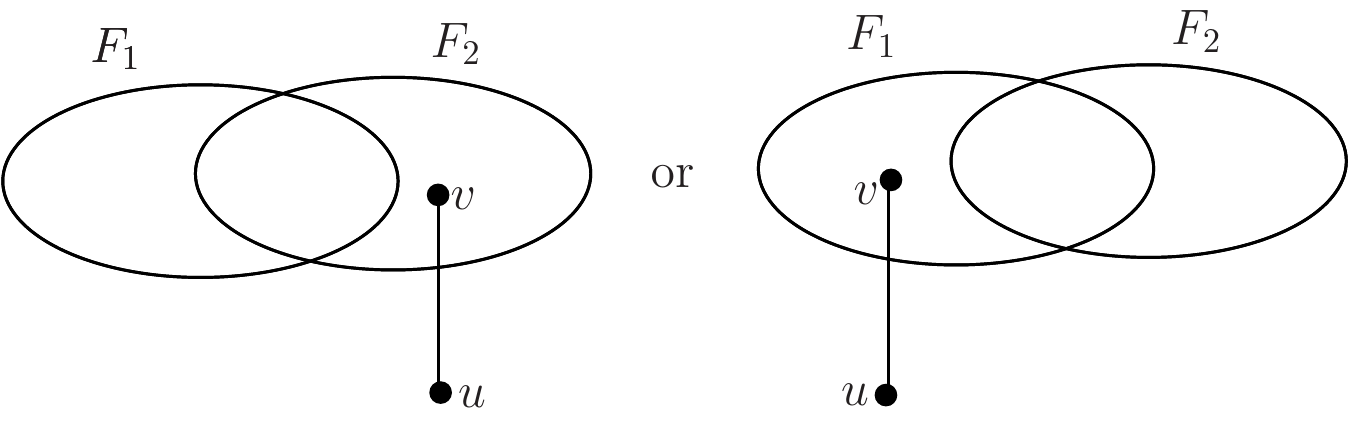}\\
  \caption{ Illustration of a distinguishable pair $(F_1, F_2)$ under the PMC model.
}\label{L111}
\end{figure}

\begin{lemma}[Dahbura and Masson \cite{DM}]\label{L01}
Let $G =(V, E)$ be a graph. For any two distinct sets $F_1, F_2 \subseteq V$, $(F_1, F_2)$ is a distinguishable pair under the PMC model if and only if there exists a vertex $u\in V-(F_1\cup F_2)$ and there exists a vertex
$v\in F_1\bigtriangleup F_2$ such that $(u, v)\in E$ (See Figure \ref{L111}).
\end{lemma}

\begin{figure}[hptb]
  \centering
  \includegraphics[width=6cm]{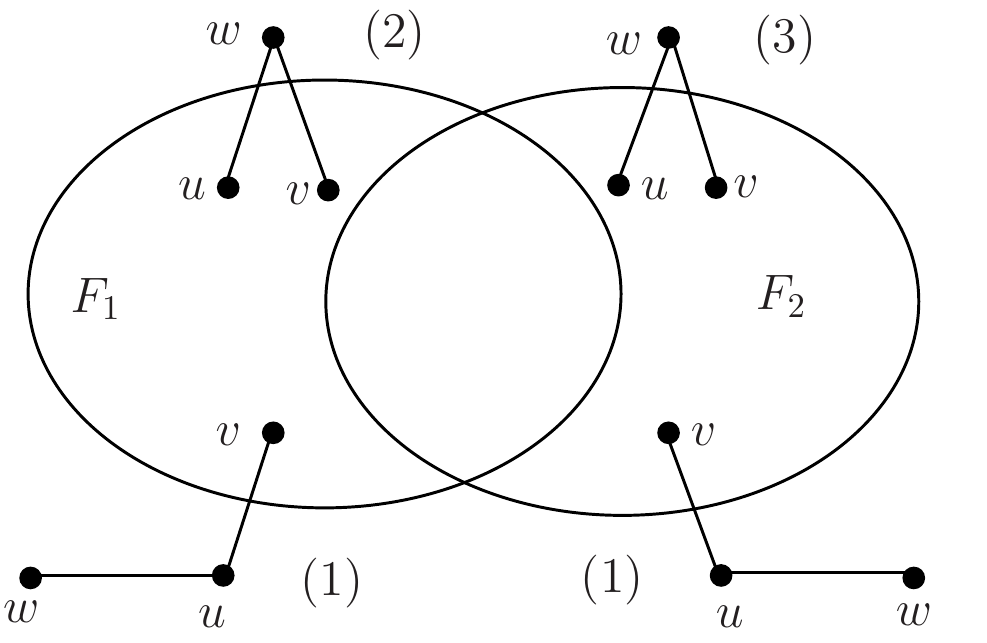}\\
  \caption{ Illustration of distinguishable sets $F_1$ and $F_2$ under the MM$^*$ model.
}\label{L22}
\end{figure}

\begin{lemma}[Sengupta and Dahbura \cite{SD}]\label{L02}
Let $G =(V, E)$ be a graph. For any two distinct sets $F_1, F_2 \subseteq V$, $F_1$ and $F_2$ are distinguishable under the MM$^*$ model if and only if any one of the following conditions is satisfied (See Figure \ref{L22}).
\begin{enumerate}
\item[{\rm (1)}] There are two vertices $u, w\in V-(F_1\cup F_2)$ and there is a vertex $v\in F_1\bigtriangleup F_2$
such that $(u, v) \in E$ and $(u, w)\in E$.

\item[{\rm (2)}] There are two vertices $u, v\in F_1-F_2$ and there is a vertex $w\in V-(F_1\cup F_2)$
such that $(u, w) \in E$ and $(v, w)\in E$.

\item[{\rm (3)}] There are two vertices $u, v\in F_2-F_1$ and there is a vertex $w\in V-(F_1\cup F_2)$
such that $(u, w) \in E$ and $(v, w)\in E$.
\end{enumerate}
\end{lemma}

Next, we introduce the diagnosability, conditional diagnosability, $g$-good-neighbor conditional diagnosability, and $R_g$-connectivity of a graph in the following statements.
\begin{definition}[Dahbura and Masson \cite{DM}]\label{D3}
For a graph $G=(V, E)$, $G$ is $t$-diagnosable if all faulty processors can be detected without replacement, provided that the number of faults does not exceed $t$. The diagnosability $t(G)$ of graph $G$ is the maximum value of $t$ such that $G$ is $t$-diagnosable.
\end{definition}

The diagnosability of multiprocessor systems, as defined above, assumes that all neighbors of any processor may fail simultaneously. However, the probability that all the neighbors of a processor fail is very small. In 2005, Lai et al. \cite{LTC} introduced conditional diagnosability under the assumption that all the neighbors of any processor in a multiprocessor system cannot be faulty at the same time.

\begin{definition}[Lai et al. \cite{LTC}]\label{D4} For a graph $G=(V, E)$, $G$ is conditionally $t$-diagnosable if $G$ is $t$-diagnosable, provided that for any processor $v \in V$, the set of faults does not contain the neighborhood $N(v)$ as a subset. The conditional diagnosability $t_c(G)$ of graph $G$ is the maximum value of $t$ such that $G$ is conditionally $t$-diagnosable.
\end{definition}

Inspired by the concept of conditional diagnosability, Peng et al. \cite{PLT} proposed $g$-good-neighbor conditional diagnosability in 2012, which extended the concept of conditional diagnosability.

\begin{definition}[Peng et al. \cite{PLT}]\label{D5} For a graph $G =(V, E)$, a faulty set $F\subseteq V$ is called a $g$-good-neighbor conditional faulty set if $|N(v) \cap (V-F)|\geq g$ for each node $v$ in $V-F$. A graph $G$ is $g$-good-neighbor conditional $t$-diagnosable if $G$ is $t$-diagnosable, provided that every faulty set is a $g$-good-neighbor conditional faulty set. The $g$-good-neighbor conditional diagnosability $t_g(G)$ of $G$ is the maximum value of $t$ such that $G$ is $g$-good-neighbor conditionally $t$-diagnosable.
\end{definition}

\begin{definition}[Yuan et al. \cite{YLM}]\label{D6}
A $g$-good-neighbor conditional cut of a graph $G$
is a $g$-good-neighbor conditional faulty set $F$ such that $G-F$ is disconnected. The minimum
cardinality of $g$-good-neighbor cuts is said to be the $R_g$-connectivity of $G$, denoted
by $\kappa^g(G)$.
\end{definition}

The same concepts as $g$-good-neighbor conditional cut and $R_g$-connectivity can be found in \cite{LGY}. We restate their definitions as follows.
\begin{itemize}
  \item Let $G$ be a connected graph. A subset $T\subset V(G)$, if any, is called an {\em $h$-vertex-cut} if $G-T$ is disconnected and has the minimum degree at least $h$. The {\em $h$-super connectivity} $\kappa^{(h)}_s(G)$ of $G$ is defined as the minimum cardinality over all $h$-vertex-cuts of $G$.
\end{itemize}
In this paper, we adopt the notation $\kappa^g(G)$.

The $R_g$-connectivity of an $(n, k)$-star graph $S_{n, k}$ was determined as follows.
\begin{lemma}[Li et al. \cite{LGY}]\label{L'5}
For $2\leq k\leq n-1$ and $n-k\leq g\leq n-2$, $\kappa^g(S_{n, k})=\dfrac{(g+1)!(n-g-1)}{(n-k)!}$.
\end{lemma}

The following lemmas give necessary and sufficient conditions for a system to be $g$-good-neighbor $t$-diagnosable under the PMC model and under the MM$^*$ model.
\begin{lemma}[Peng et al. \cite{PLT} and Dahbura and Masson \cite{DM}]\label{L1}
A graph $G =(V, E)$ is $g$-good-neighbor $t$-diagnosable under the
PMC model if and only if there is an edge $(u, v)\in E$ with $u\in V-(F_1\cup F_2)$ and
$v\in F_1\bigtriangleup F_2$ for each distinct pair of $g$-good-neighbor conditional faulty sets $F_1$ and $F_2$ of $V$ with $|F_1|\leq t$ and $|F_2|\leq t$ (See Figure \ref{L111}).
\end{lemma}

\begin{lemma}[Sengupta and Dahbura \cite{SD} and Yuan et al. \cite{YLM}]\label{L2}
A graph $G =(V, E)$ is $g$-good-neighbor $t$-diagnosable under the
MM$^*$ model if and only if each distinct pair of $g$-good-neighbor conditional faulty sets $F_1$ and
$F_2$ of $V$ with $|F_1|\leq t$ and $|F_2|\leq t$ satisfies one of the following conditions (See Figure \ref{L22}).
\begin{enumerate}
\item[{\rm (1)}] There are two vertices $u, w\in V-(F_1\cup F_2)$ and there is a vertex $v\in F_1\bigtriangleup F_2$
such that $(u, v) \in E$ and $(u, w)\in E$.

\item[{\rm (2)}] There are two vertices $u, v\in F_1-F_2$ and there is a vertex $w\in V-(F_1\cup F_2)$
such that $(u, w) \in E$ and $(v, w)\in E$.

\item[{\rm (3)}] There are two vertices $u, v\in F_2-F_1$ and there is a vertex $w\in V-(F_1\cup F_2)$
such that $(u, w) \in E$ and $(v, w)\in E$.
\end{enumerate}
\end{lemma}

Next, we introduce the split-graph, proposed in \cite{LGY},  which will be used in the proof of our main results.

\begin{figure}[hptb]
  \centering
  \includegraphics[width=10cm]{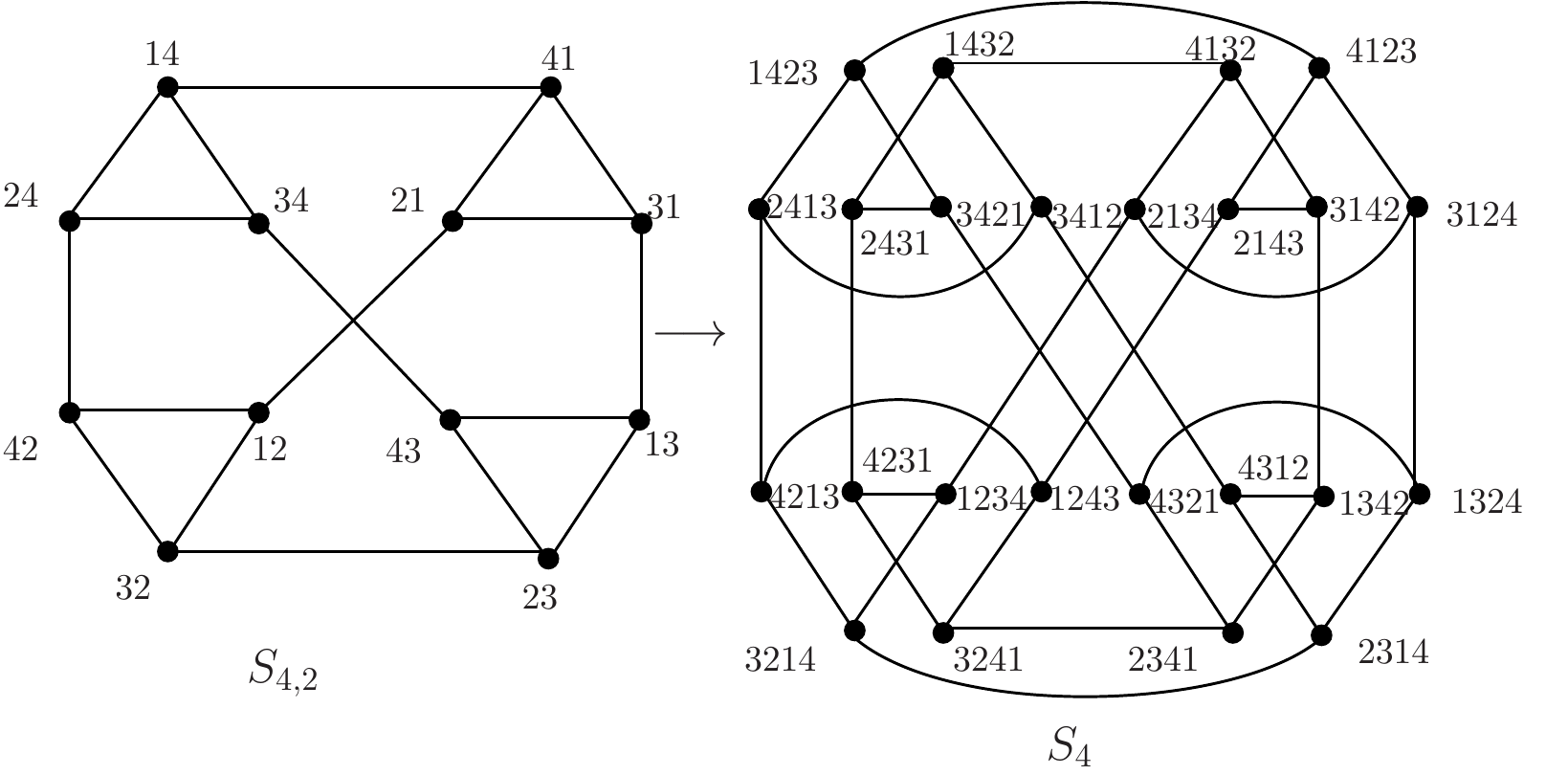}\\
  \caption{A $(4,2)$-star graph $S_{4,2}$ and its $2$-split graph $S_{4,2}^2$, which is isomorphic to a star $S_4$.
}\label{D003}
\end{figure}

\begin{definition}[Li et al. \cite{LGY}]\label{D03}
Let $G$ be a graph and $t$ be a positive integer. A $t$-split graph $G^t$ of $G$ is
a graph obtained from $G$ by replacing each vertex $x$ by a set $V_x$ of $t$ independent vertices,
and replacing each edge $e = xy$ by a perfect matching $E_e$ between $V_x$ and $V_y$ (See Figure \ref{D003}).
\end{definition}

Li et al. \cite{LGY} obtained the relationship between $S_n$ and $S_{n, k}$ as follows.

\begin{lemma}[Li et al. \cite{LGY}]\label{L'3}
For any $k$ with $2\leq k\leq n-1$, there is an $(n-k)!$-split graph of $S_{n, k}$ that is isomorphic to a star graph $S_n$.
\end{lemma}

The following lemma is very useful in proving our main results.

\begin{lemma}[Li and Lu \cite{LL}]\label{L11}
Let $H$ be a subgraph of $S_n$ for $n\geq2$ and $\delta(H)\geq g$, where $g\in I_{n-1}$. Then $|V(H)|\geq (g+1)!$.
\end{lemma}

\section{Main Results}\label{4}
\begin{lemma}\label{L11'}
Let $H$ be a subgraph of $S_{n, k}$ for $2\leq k\leq n-1$. If $\delta(H)\geq g$, where $g\in I_{n-1}$, then $|V(H)|\geq \dfrac{(g+1)!}{(n-k)!}$.
\end{lemma}
\proof
Suppose that $H$ is a subgraph of $S_{n, k}$. Since $2\leq k\leq n-1$, by Lemma \ref{L'3}, there is an $(n-k)!$-split graph of $S_{n, k}$, denoted by $S_{n, k}^{(n-k)!}$ which is isomorphic to a star graph $S_n$. Let $H^{(n-k)!}$ be an $(n-k)!$-split graph of $H$ satisfying that it is also a subgraph of $S_{n, k}^{(n-k)!}$. Thus, $|V(H^{(n-k)!})|=(n-k)!|V(H)|$. Since $\delta(H)\geq g$, we have $\delta( H^{(n-k)!})\geq g$. By Lemma \ref{L11}, $|V(H^{(n-k)!})|\geq (g+1)!$. Therefore, we have $|V(H)|\geq \dfrac{(g+1)!}{(n-k)!}$.

We obtain the desired result. $\qed$

Note that $S_{n, k}$ is $(n-1)$-regular. Then $t_{n-1}(S_{n, k})=0$ under the PMC model and the MM$^*$ model. In the following, we assume that $g\leq n-2$.

Now, we consider the upper bound of the $g$-good-neighbor conditional diagnosability of $(n, k)$-star network $S_{n, k}$ for $2\leq k\leq n-1$ and $n-k\leq g\leq n-2$ under the PMC model and the MM$^*$ model.

\begin{figure}[hptb]
  \centering
  \includegraphics[width=6cm]{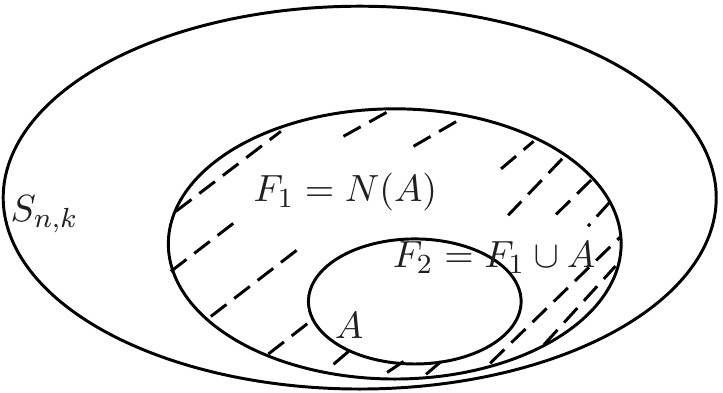}\\
  \caption{ Illustration of $F_1$ and $F_2$.
}\label{P11}
\end{figure}

\begin{lemma}\label{p1}
For $2\leq k\leq n-1$ and $n-k\leq g\leq n-2$, we have $t_g(S_{n, k})\leq\dfrac{(g+1)!(n-g)}{(n-k)!}-1$ under the PMC model and the MM$^*$ model.
\end{lemma}
\proof
Let $$A=\{p_1p_2\ldots p_{k-(n-g-1)}12\ldots(n-g-1)\in V(S_{n, k})\mid p_i\in I_n-I_{n-g-1},\ i\in I_{k-(n-g-1)}\}, $$ $F_1=N(A)$ and $F_2=F_1\cup A$ (see Figure \ref{P11}).
Then $|A|=\dfrac{(g+1)!}{(n-k)!}$. For $F_1=N(A)$, we have
\begin{equation*}
  F_1=\{ip_2\ldots p_{k-(n-g-1)}1\ldots(i-1)p_1(i+1)\ldots(n-g-1)\in V(S_{n, k})|p_j\in I_n-I_{n-g-1},\ j\in I_{n-g-1}\}.
\end{equation*}
Choose a vertex $u=p_1p_2\ldots p_{k-(n-g-1)}12\ldots(n-g-1)\in A$. By Definition \ref{D02}, we have
\begin{eqnarray*}
N(u)\cap A&=&\{p_ip_2\ldots p_{i-1}p_1p_{i+1}\ldots p_{k-(n-g-1)}12\ldots(n-g-1)\mid i\in I'_{k-(n-g-1)}\}\cup\\
&&\{qp_2\ldots p_{k-(n-g-1)}12\ldots(n-g-1)\mid q\in I_n-I_{n-g-1}-\{p_i\mid i\in I_{k-(n-g-1)}\}\},
\end{eqnarray*}
$$N(u)\cap F_1=\{ip_2\ldots p_{k-(n-g-1)}1\ldots(i-1)p_1(i+1)\ldots(n-g-1)\in V(S_{n, k}) \mid i\in I_{n-g-1}\}. $$
Then $|N(u)\cap A|=g$ and $|N(u)\cap F_1|=n-g-1$.
By the definition of $S_{n, k}$, no two vertices in $A$ share a common neighbor in $F_1$. It follows that $|F_1|=|A|(n-g-1)=\dfrac{(g+1)!(n-g-1)}{(n-k)!}$ and $|F_2|=|F_1|+|A|=\dfrac{(g+1)!(n-g)}{(n-k)!}$.
Note that $A=F_1\bigtriangleup F_2$ and $F_1=N(A)$. Since there is no edge between $F_1\bigtriangleup F_2$ and $V(S_{n, k})-(F_1\cup F_2)$, by Lemmas \ref{L01} and \ref{L02}, we conclude that $F_1$ and $F_2$ are indistinguishable under the PMC model and the MM$^*$ model.

Now, we verify that both $F_1$ and $F_2$ are $g$-good-neighbor conditional faulty sets.

Suppose $u\in V(S_{n, k})-F_2$.
If $u \notin N(F_2)$, then $N(u)\subseteq V(S_{n, k})-F_2$ and $|N(u)|=n-1>g$.
If $u \in N(F_2)$, then there exists $v\in F_1$ such that $uv\in E(S_{n, k})$. We can assume that $v=ip_2\ldots p_{k-(n-g-1)}1\ldots (i-1)p_1(i+1)\ldots (n-g-1)$ for some fixed $i\in I_{n-g-1}$. Thus, $u$ has three forms, i.e.,
\begin{equation*}\label{1111}
 jp_2\ldots p_{k-(n-g-1)}1\ldots (j-1)i(j+1)\ldots (i-1)p_1(i+1)\ldots (n-g-1),\quad j\in I_{n-g-1},\ i\neq j;
\end{equation*}
\begin{equation*}\label{2222}
  p_jp_2\ldots p_{j-1}ip_{j+1}\ldots p_{k-(n-g-1)}1\ldots (i-1)p_1(i+1)\ldots (n-g-1),\quad j\in I'_{k-(n-g-1)};
\end{equation*}
\begin{equation*}\label{3333}
  qp_2\ldots p_{k-(n-g-1)}1\ldots (i-1)p_1(i+1)\ldots (n-g-1), q\in I_n-I_{n-g-1}-\{p_1, p_2, \ldots, p_{k-(n-g-1)}\}.
\end{equation*}
No matter which form $u$ has, $N(u)\cap F_1=\{v\}$. Since $S_{n, k}$ is $(n-1)$-regular, $|N(u)-F_2|=(n-1)-1\geq g$.
Thus, $F_2$ is a $g$-good-neighbor conditional faulty set.

Suppose $u\in V(S_{n, k})-F_1$.
If $u \in V(LTQ_n)-F_2$, then we obtain the desired result by the same proof as above.
If $u \notin V(LTQ_n)-F_2$, then $u \in A$. Note that $S_{n, k}[A]$ is $g$-regular. Thus, $|N(u)\cap A|=g$. So $F_1$ is a $g$-good-neighbor conditional faulty set.

By Lemmas \ref{L1} and \ref{L2}, we have $t_g(S_{n, k})\leq\dfrac{(g+1)!(n-g)}{(n-k)!}-1$ under the PMC model and the MM$^*$ model. $\qed$

Next, we consider the lower bound of the $g$-good-neighbor conditional diagnosability of $(n, k)$-star network $S_{n, k}$  for $2\leq k\leq n-1$ and $n-k\leq g\leq n-2$ under the PMC model and the MM$^*$ model, respectively.

\begin{lemma}\label{p2}
For $2\leq k\leq n-1$ and $n-k\leq g\leq n-2$ with $n\geq4$,  we have $t_g(S_{n, k})\geq\dfrac{(g+1)!(n-g)}{(n-k)!}-1$ under the PMC model.
\end{lemma}
\proof
Suppose that $F_1$ and $F_2$ are any two distinct $g$-good-neighbor conditional faulty sets and they are indistinguishable. We will prove the lemma by showing that $|F_1|\geq \dfrac{(g+1)!(n-g)}{(n-k)!}$ or $|F_2|\geq \dfrac{(g+1)!(n-g)}{(n-k)!}$.

If $V(S_{n, k})=F_1\cup F_2$, then $|F_1|\geq \dfrac{1}{2}|V(S_{n, k})|\geq\dfrac{n!}{2(n-k)!}\geq \dfrac{(g+1)!(n-g)}{(n-k)!}$ or $|F_2|\geq \dfrac{1}{2}|V(S_{n, k})|\geq\dfrac{n!}{2(n-k)!}\geq \dfrac{(g+1)!(n-g)}{(n-k)!}$ for $n\geq 4$.

Now, we suppose $V(S_{n, k})\neq F_1\cup F_2$. Since $F_1$ and $F_2$ are indistinguishable, there are no edges between $V(S_{n, k})-(F_1\cup F_2)$ and $F_1\bigtriangleup F_2$ by Lemma \ref{L01}. Thus, $S_{n, k}-(F_1\cap F_2)$ is disconnected. Without loss of generality, we assume that $F_2-F_1\neq \emptyset$. Note that $F_1$ and $F_2$ are both $g$-good-neighbor conditional faulty sets. Then $\delta(S_{n, k}[V(S_{n, k})-(F_1\cup F_2)])\geq g$ and $\delta(S_{n, k}[F_2-F_1])\geq g$. Thus, we have $|F_2-F_1|\geq \dfrac{(g+1)!}{(n-k)!}$ by Lemma \ref{L11'}. If $F_1-F_2\neq\emptyset$, then $\delta(S_{n, k}[F_1-F_2])\geq g$. Therefore, $F_1\cap F_2$ is a $g$-good-neighbor conditional cut of $S_{n, k}$. By Lemma \ref{L'5}, we obtain that $|F_1\cap F_2|\geq \dfrac{(g+1)!(n-g-1)}{(n-k)!}$. Hence, $|F_2|=|F_2-F_1|+|F_1\cap F_2|\geq \dfrac{(g+1)!}{(n-k)!}+\dfrac{(g+1)!(n-g-1)}{(n-k)!}=\dfrac{(g+1)!(n-g)}{(n-k)!}$.

This completes the proof of Lemma \ref{p2}. $\qed$

\begin{lemma}\label{p3}
For $2\leq k\leq n-1$ and $n-k\leq g\leq n-2$ with $n\geq4$, we have $t_g(S_{n, k})\geq\dfrac{(g+1)!(n-g)}{(n-k)!}-1$ under the MM$^*$ model.
\end{lemma}
\proof
If $g=1$, then $k=n-1$. Note that $S_{n, n-1}$ is isomorphic to $S_n$. By Theorems \ref{T1} and \ref{T2}, $t_1(S_{n, n-1})=2n-3$ for $n\geq4$. The result holds.

Now assume $g\geq 2$. Suppose that $F_1$ and $F_2$ are any two distinct $g$-good-neighbor conditional faulty sets and they are indistinguishable. We will prove the lemma by showing that $|F_1|\geq \dfrac{(g+1)!(n-g)}{(n-k)!}$ or $|F_2|\geq \dfrac{(g+1)!(n-g)}{(n-k)!}$.

If $V(S_{n, k})=F_1\cup F_2$, then $|F_1|\geq \dfrac{1}{2}|V(S_{n, k})|\geq\dfrac{n!}{2(n-k)!}\geq \dfrac{(g+1)!(n-g)}{(n-k)!}$ or $|F_2|\geq \dfrac{1}{2}|V(S_{n, k})|\geq\dfrac{n!}{2(n-k)!}\geq \dfrac{(g+1)!(n-g)}{(n-k)!}$ for $n\geq 4$.

Now, we suppose $V(S_{n, k})\neq F_1\cup F_2$. Without loss of generality, we assume that $F_2-F_1\neq \emptyset$. We shall show that there is no edge between $F_1\bigtriangleup F_2$ and $V(S_{n, k})-(F_1\cup F_2)$. Otherwise, there exists an edge $uv\in E(S_{n, k})$, where $u\in F_1\bigtriangleup F_2$ and $v\in V(S_{n, k})-(F_1\cup F_2)$. Without loss of generality, we can assume that $u\in F_2-F_1$. Since $F_1$ is a $g$-good-neighbor conditional faulty set with $g\geq 2$, $v$ has at least two neighbors in $S_{n, k}-F_1$. Thus, $v$ has a neighbor $w~(w\neq u)$ in $F_2-F_1$ or $V(S_{n, k})-(F_1\cup F_2)$, which contradicts Lemma \ref{L02}. By the same discussion as Lemma \ref{p2}, we complete the proof of this lemma. $\qed$

Now, we get the result below by Lemmas \ref{p1}, \ref{p2}, and \ref{p3} directly.

\begin{thm}\label{main1}
Let $S_{n, k}$ be $(n, k)$-star networks with $2\leq k\leq n-1$ and $n\geq4$. Then the $g$-good-neighbor conditional diagnosability of $S_{n, k}$ under the PMC model and the MM$^*$ model are both $t_g(S_{n, k})=\dfrac{(g+1)!(n-g)}{(n-k)!}-1$ for $n-k\leq g\leq n-2$.
\end{thm}

Now we discuss the other unknown cases in Table \ref{tab:1}.
\begin{thm}\label{main2}
The $g$-good-neighbor conditional diagnosability of $S_{n, 1}$ with $n\geq4$ under the PMC model and the MM$^*$ model are both $t_g(S_{n, 1})=n-g-1$ for $\left\lfloor\dfrac{n}{2}\right\rfloor\leq g\leq n-2$. What is more, the $1$-good-neighbor conditional diagnosabilities of $S_{3, 1}$ under the PMC model and the MM$^*$ model are $1$ and $0$, respectively.
\end{thm}
\proof
When $k=1$, $S_{n, 1}$ is isomorphic to $K_n$ which is a complete graph.

First, we consider the $1$-good-neighbor conditional diagnosability of $S_{3, 1}$ under the PMC model and the MM$^*$ model.

Note that any $1$-good-neighbor conditional faulty set of $S_{3, 1}$ contains at most one vertex. We have $t_1(S_{3, 1})\leq1$. Let $F_1$ and $F_2$ are any two distinct $1$-good-neighbor conditional faulty sets of $S_{3, 1}$ with $|F_1|\leq1$ and $|F_2|\leq1$. Since $V(S_{3, 1})-(F_1\cup F_2)\neq\emptyset$, there exists one edge between $F_1\bigtriangleup F_2$ and $V(S_{n, k})-(F_1\cup F_2)$. By Lemma \ref{L01}, we know $F_1$ and $F_2$ are distinguishable under the PMC model. Thus, $t_1(S_{3, 1})=1$ under the PMC model.

Suppose $V(S_{3, 1})=\{u_1, u_2, u_3\}$. Let $F_1=\{u_1\}$ and $F_2=\{u_2\}$. Then $F_1\bigtriangleup F_2=\{u_1, u_2\}$ and $V(S_{3, 1})-(F_1\cup F_2)=\{u_3\}$. Obviously, $F_1$ and $F_2$ are $1$-good-neighbor conditional faulty sets of $S_{3, 1}$. By Lemma \ref{L02}, $F_1$ and $F_2$ are indistinguishable under the MM$^*$ model. Hence, $t_1(S_{3, 1})=0$ under the MM$^*$ model.

Next, we consider the $g$-good-neighbor conditional diagnosability of $S_{n, 1}$ with $n\geq4$ for $\left\lfloor\dfrac{n}{2}\right\rfloor\leq g\leq n-2$ under the PMC model and the MM$^*$ model.

Note that the $g$-good-neighbor conditional faulty set of $S_{n, 1}$ contains at most $n-g-1$ vertices. We have $t_g(S_{n, 1})\leq n-g-1$ under the PMC model and the MM$^*$ model.

To proceed, we show that $t_g(S_{n, 1})\geq n-g-1$ under the PMC model and the MM$^*$ model.

Suppose that $F_1$ and $F_2$ are any two distinct $g$-good-neighbor conditional faulty sets and $|F_i|\leq n-g-1$ for each $i\in I_2$. We will show that they are distinguishable. Without loss of generality, we assume that $F_2-F_1\neq \emptyset$.

If $V(S_{n, 1})=F_1\cup F_2$, then $n=|V(S_{n, 1})|=|F_1\cup F_2|=|F_1|+|F_2|-|F_1\cap F_2|\leq |F_1|+|F_2|\leq 2(n-g-1)\leq 2(n-\left\lfloor\dfrac{n}{2}\right\rfloor-1)<n$, which is a contradiction.

Now, we suppose $V(S_{n, 1})\neq F_1\cup F_2$. We have
\begin{equation}\label{eee}
\begin{split}
 |V(S_{n, 1})|-|F_1\cup F_2|&= |V(S_{n, 1})|-|F_1|-|F_2|+|F_1\cap F_2|\\
&\geq n-2(n-g-1)+|F_1\cap F_2|\\
&= 2(g+1)-n+|F_1\cap F_2|\\
&\geq 2(\left\lfloor\dfrac{n}{2}\right\rfloor+1)-n+|F_1\cap F_2|\\
&\geq 1+|F_1\cap F_2|.
\end{split}
\end{equation}
Note that $S_{n, 1}$ is isomorphic to $K_n$. Under the PMC model, by Lemma \ref{L01}, we know that $F_1$ and $F_2$ are distinguishable. Thus, we have $t_g(S_{n, 1})\geq n-g-1$ under the PMC model.

\begin{figure}[hptb]
  \centering
  \includegraphics[width=10cm]{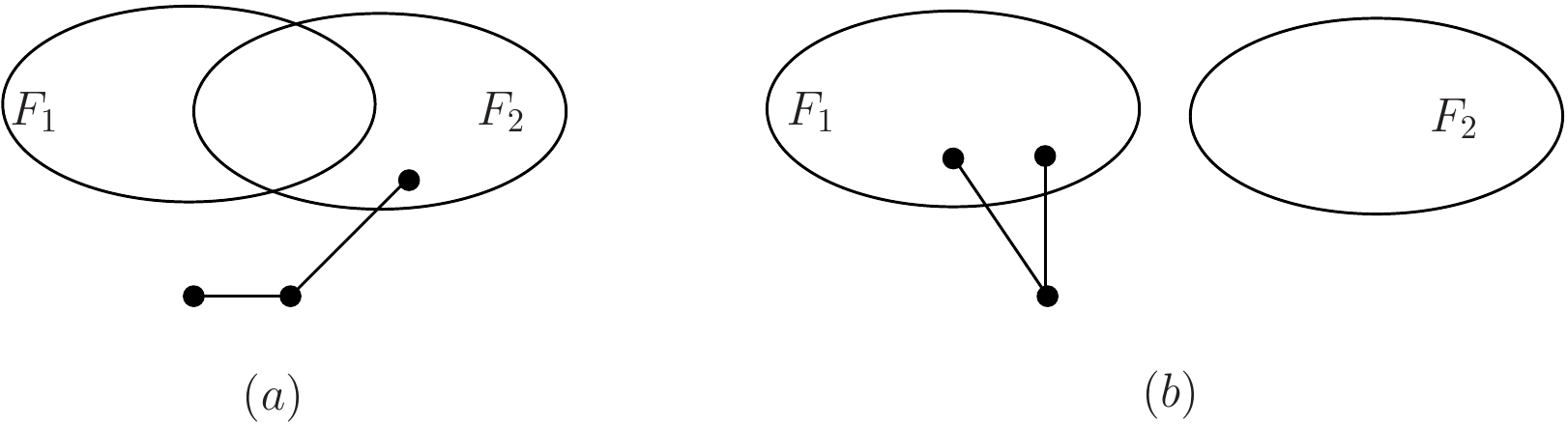}\\
  \caption{ Illustration of the proof of Theorem \ref{main2}.
}\label{1ab}
\end{figure}

If $|V(S_{n, 1})-(F_1\cup F_2)|\geq 2$, then by Lemma \ref{L02}(1), we know that $F_1$ and $F_2$ are distinguishable. See Figure \ref{1ab}(a).

If $|V(S_{n, 1})-(F_1\cup F_2)|= 1$, then by (\ref{eee}), we have $n$ is odd and $F_1\cap F_2=\emptyset$. We can also know that $|F_1|=|F_2|=\left\lfloor\dfrac{n}{2}\right\rfloor\geq 2$. Thus, by Lemma \ref{L02}(3), we know that $F_1$ and $F_2$ are distinguishable. See Figure \ref{1ab}(b).

Thus, we have $t_g(S_{n, 1})\geq n-g-1$ under the MM$^*$ model.

As mentioned above, we obtain the desired result. $\qed$

Finally, we discuss $t_1(S_{n, 2})$ under the MM$^*$ model.
\begin{thm}\label{main3}
Under the MM$^*$ model, we have
\begin{equation*}\label{ttt}
  t_1(S_{n, 2})=
\left\{
  \begin{array}{lll}
   n-1 & \hbox{if $n\geq4$}; \\
    \\
    1 & \hbox{if $n=3$}.
  \end{array}
\right.
\end{equation*}
\end{thm}
\proof
First, we consider $t_1(S_{3, 2})$ under the MM$^*$ model. Note that $S_{3, 2}$ is a cycle with six vertices. Suppose $V(S_{3, 2})=\{u_1, u_2, \ldots, u_6\}$ and $E(S_{3, 2})=\{u_1u_2, u_2u_3, u_3u_4, u_4u_5, u_5u_6, u_1u_6\}$.

\begin{figure}[hptb]
  \centering
  \includegraphics[width=6cm]{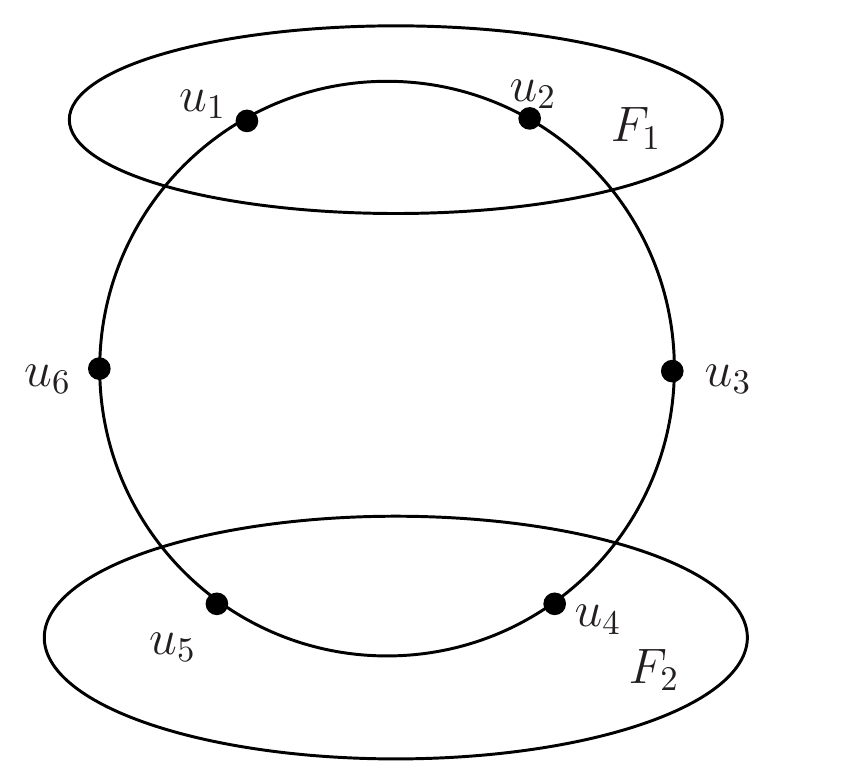}\\
  \caption{ Illustration of the proof of $t_1(S_{3, 2})\leq 1$ under the MM$^*$ model.
}\label{c6s}
\end{figure}

Let $F_1=\{u_1, u_2\}$ and $F_2=\{u_4, u_5\}$. Obviously, $F_1$ and $F_2$ are $1$-good-neighbor conditional faulty sets of $S_{3, 2}$. By Lemma \ref{L02}, $F_1$ and $F_2$ are indistinguishable under the MM$^*$ model. See Figure \ref{c6s}. Thus, we have $t_1(S_{3, 2})\leq 1$ under the MM$^*$ model.

On the other hand, by transitivity of $S_{3, 2}$, we may suppose that $F_1=\{u_1\}$ and $F_2=\{u_j\}$ (or $F_2=\emptyset$) are any two distinct $1$-good-neighbor conditional faulty sets of $S_{3, 2}$.

When $F_2=\{u_j\}$, we assume that $2\leq j\leq6$.

\begin{figure}[hptb]
  \centering
  \includegraphics[width=10cm]{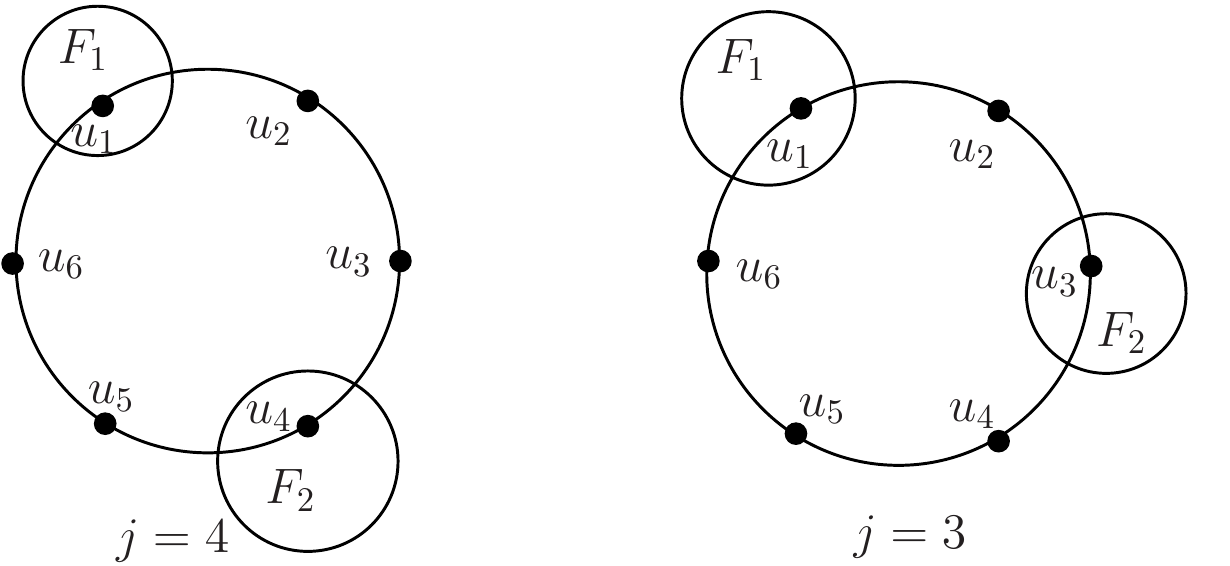}\\
  \caption{ Illustration of the proof of $t_1(S_{3, 2})\geq 1$ under the MM$^*$ model.
}\label{c6abc}
\end{figure}

If $4\leq j\leq6$, then $u_2, u_3\in V(S_{3, 2})-(F_1\cup F_2)$. Since $u_1\in F_1\bigtriangleup F_2$ and $u_1u_2, u_2u_3\in E(S_{3, 2})$, we have that $F_1$ and $F_2$ are distinguishable under the MM$^*$ model. See Figure \ref{c6abc} for $j=4$.

If $2\leq j\leq3$, then $u_5, u_6\in V(S_{3, 2})-(F_1\cup F_2)$. Since $u_1\in F_1\bigtriangleup F_2$ and $u_1u_6, u_6u_5\in E(S_{3, 2})$, we have that $F_1$ and $F_2$ are distinguishable under the MM$^*$ model.
See Figure \ref{c6abc} for $j=3$.

When $F_2=\emptyset$, the result also holds.

By Lemma \ref{L2}, we have $t_1(S_{3, 2})\geq 1$ under the MM$^*$ model.
Therefore, $t_1(S_{3, 2})= 1$ under the MM$^*$ model.

Now, we consider $t_1(S_{n, 2})$ under the MM$^*$ model, where $n\geq4$. Our discussion is divided into two steps.

\medskip

\noindent \textbf{Step 1:} Show that $t_1(S_{n, 2})\geq n-1$ under the MM$^*$ model, where $n\geq4$.

Suppose that $F_1$ and $F_2$ are any two distinct $1$-good-neighbor conditional faulty sets and they are indistinguishable. We will show that $|F_1|\geq n$ or $|F_2|\geq n$.

If $V(S_{n, 2})=F_1\cup F_2$, then $|F_1|\geq \dfrac{1}{2}|V(S_{n, 2})|=\dfrac{n(n-1)}{2}\geq n$ or $|F_2|\geq \dfrac{1}{2}|V(S_{n, 2})|=\dfrac{n(n-1)}{2}\geq n$ for $n\geq4$.

Now, we suppose $V(S_{n, 2})\neq F_1\cup F_2$. Without loss of generality, we assume that $F_2-F_1\neq \emptyset$.

If $|F_1\cap F_2|\geq n-1$, then $|F_2|=|F_2-F_1|+|F_1\cap F_2|\geq 1 + (n-1)=n$.

Now, we suppose $|F_1\cap F_2|\leq n-2$. Let $W_1, \ldots, W_c$ be the components of $S_{n, 2}-(F_1\cup F_2)$ such that $|V(W_1)|\leq \ldots\leq |V(W_c)|$, where $c\geq 1$. For any component $W_i$, if $|V(W_i)|\geq 2$, then there is no edge between $W_i$ and $F_1\bigtriangleup F_2$. Otherwise, it contradicts the fact that $F_1$ and $F_2$ are indistinguishable by Lemma \ref{L02}.

If $|V(W_c)|\geq 2$, then $F_1\cap F_2$ is a cut of $S_{n, 2}$.  Since the connectivity of $S_{n, 2}$ is $n-1$, we have $|F_1\cap F_2|\geq n-1$, which contradicts the assumption that $|F_1\cap F_2|\leq n-2$.

Next, we assume that $|V(W_c)|=1$. If $F_1-F_2=\emptyset$, then $S_{n, 2}-F_1-F_2=S_{n, 2}-F_2$. The vertex in $W_1$ is one isolated vertex in $S_{n, k}-F_2$, which contradicts the fact that $F_2$ is a $1$-good neighbor conditional faulty set. We suppose $F_1-F_2\neq\emptyset$. Let $W=\bigcup_{1\leq i\leq c}V(W_i)$. Then $V(S_{n, 2})=W\cup(F_1\cup F_2)$. Arbitrarily choose a vertex $w\in W$. Then, $N(w)\subseteq F_1\cup F_2$. Since $F_1$ and $F_2$ are indistinguishable, $|N(w)\cap(F_2-F_1)|\leq 1$ and $|N(w)\cap(F_1-F_2)|\leq 1$ by Lemma \ref{L02}.
Owing to the fact that $F_1$ and $F_2$ are $1$-good-neighbor conditional faulty sets, we have $|N(w)\cap(F_2-F_1)|=|N(w)\cap(F_1-F_2)|=1$ and $|N(w)\cap(F_1\cap F_2)|=n-3\leq |F_1\cap F_2|$. Thus,
\begin{eqnarray*}
\sum _{w\in W}|N(w)\cap(F_1\cap F_2)|&=&|W|(n-3)\\
&\leq& \sum _{x\in F_1\cap F_2}d(x)\\
&\leq& |F_1\cap F_2|(n-1)\\
&\leq& (n-2)(n-1).
\end{eqnarray*}
It follows that $|W|\leq \dfrac{(n-2)(n-1)}{n-3}=n+\dfrac{2}{n-3}\leq n+2$ when $n\geq 4$. Thus,
\begin{eqnarray*}
|F_1|+|F_2|&=& |V(S_{n, 2})|+|F_1\cap F_2|-|W|\\
&\geq& n(n-1)+(n-3)-(n+2)\\
&=& n^2-n-5.
\end{eqnarray*}
Therefore, for $n\geq4$, we have
\begin{eqnarray*}
\max\{|F_1|, |F_2|\}&\geq& \left\lceil\dfrac{|F_1|+|F_2|}{2}\right\rceil\\
&\geq& \left\lceil\dfrac{n^2-n-5}{2}\right\rceil\\
&\geq& n.
\end{eqnarray*}

Based on the above discussion, we obtain the desired result.

\medskip

\noindent \textbf{Step 2:} Show that $t_1(S_{n, 2})\leq n-1$ under the MM$^*$ model, where $n\geq4$.

\begin{figure}[hptb]
  \centering
  \includegraphics[width=8cm]{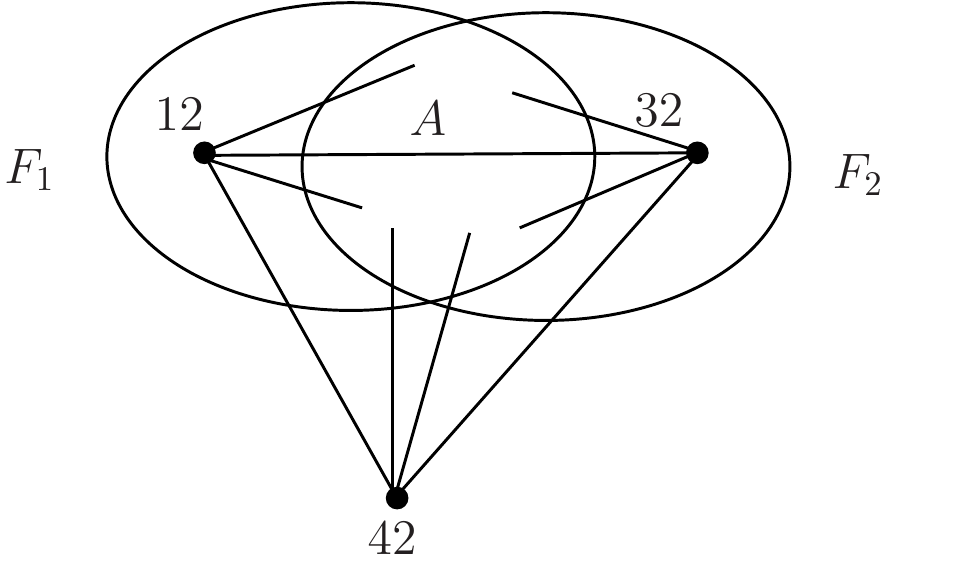}\\
  \caption{ Illustration of $F_1$ and $F_2$.
}\label{step2}
\end{figure}

By using the construction of Chang et al. \cite{CDH}, let $$A=N(\{12, 32, 42\}), $$ $$F_1=A\cup \{12\}$$  and
$$
F_2=A\cup \{32\}.
$$
Then $|A|=(n-4)+3=n-1$ and $|F_1|=|F_2|=n$. See Figure \ref{step2}.
We conclude that $F_1$ and $F_2$ are indistinguishable $1$-good-neighbor conditional faulty sets under the MM$^*$ model from the proof of \cite{CDH}.
By Lemma \ref{L2}, we have $t_1(S_{n, 2})\leq n-1$ under the MM$^*$ model, where $n\geq4$.

\medskip

By Steps 1 and 2, we conclude that $t_1(S_{n, 2})= n-1$ under the MM$^*$ model, where $n\geq4$.

This completes the proof of Theorem \ref{main3}. $\qed$

\section{Conclusions}\label{6}
In this paper, we determined the $g$-good-neighbor conditional diagnosability of $(n, k)$-star networks $S_{n, k}$ for all the remaining cases with $1\leq k\leq n-1$ and $n\geq3$ for $1\leq g\leq n-k$ under the PMC model and the MM$^*$ model (see Table \ref{tab:2}). Future research on this topic will involve studying the $g$-good-neighbor conditional diagnosability of many network topologies.
\begin{table}[htbp]\tiny
\caption{\quad \label{tab:2} The $g$-good-neighbor conditional diagnosability of $S_{n, k}$ under the PMC model and the MM$^*$ model}
\begin{center}
\extrarowheight=4pt
\renewcommand{\arraystretch}{1.3}
\begin{tabular}{|c|c|c|c|c|}\hline
 &$k=1$ &$k=2$ &$3\leq k\leq n-2$& $k=n-1$ \\
\hline
&$\left\lceil\dfrac{n}{2}\right\rceil-1 ~(n\geq4)$ \cite{XL} & $n$ (PMC) \cite{XL} &&\\
$g=1$ &$0$ ($n=3$, {\rm MM}$^*$)~(Theorem \ref{main2}) & $n-1$ ($n\geq4$, {\rm MM}$^*$)~(Theorem \ref{main3}) &{$n+k-2$ \cite{XL}}&$2n-3$~\cite{LL,XL}\\
&$1$ ($n=3$, PMC)~(Theorem \ref{main2}) &$1$ ($n=3$, {\rm MM}$^*$)~(Theorem \ref{main3})&&\\
\hline
$2\leq g\leq \left\lfloor\dfrac{n}{2}\right\rfloor-1$ &$\left\lceil\dfrac{n}{2}\right\rceil-1$ \cite{XL} &\multicolumn{2}{|c|}{{\multirow{2}{*}{$n+g(k-1)-1$ \cite{XL}}}}&{\multirow{3}{*}{}}\\
\cline{1-2}
$\left\lfloor\dfrac{n}{2}\right\rfloor\leq g\leq n-k$ &$n-g-1$~(Theorem \ref{main2}) &\multicolumn{2}{|c|}{}&{$(n-g)(g+1)!-1$~\cite{LL}}\\\cline{1-4}

{$n-k\leq g\leq n-2$} &Nonexistence &\multicolumn{2}{|c|}{$\dfrac{(g+1)!(n-g)}{(n-k)!}-1~ (n\geq4)$~(Theorem \ref{main1})}&\\
\hline
{ $g=n-1$} &\multicolumn{4}{|c|}{$0$}\\
\hline
\end{tabular}
\end{center}
\end{table}

\end{CJK*}

\end{document}